\documentstyle{article}
\newtheorem{theorem}{THEOREM}[section]

\newtheorem{corollary}[theorem]{Corollary}
\newtheorem{lemma}[theorem]{Lemma}

\newtheorem{proposition}[theorem]{Proposition}
\def\<{\prec}

\def\H{{\cal H}}

\def\Omegabar{\overline{\Omega}}

\def\CC{{\rm\kern.24em\vrule
width.02em height1.4ex
depth-.05ex\kern-.26em C}}
\def\QQ{{\rm\kern.24em\vrule width.02em
height1.4ex depth-.05ex\kern-.26em Q}}
\def\RR{{\rm I\kern-.2em R}}
\def\HH{{\rm I\kern-.2em H}}
\def\ZZ{{\rm\kern.26em\vrule width.02em
height0.5ex depth0ex\kern.04em\vrule width.02em
height1.47ex depth-1ex\kern-.34em Z}}
\def\Ibb#1{{\rm I\kern-.23em#1}}
\def\Ib#1{{\rm I\kern-.25em#1}}
\def\k#1{\kern#1em}
\def\vb#1{\vrule width.02em height1.4ex depth-.05ex}

\def\11{{\rm\k{.45}\vb0\k{-.142}1}}
\def\c{\rfloor}

\def\epf{\hskip.2in\vrule width.4pt height6.65pt
depth.15pt\vrule
width2.5pt height6.65pt depth-6.25pt\hskip-2.5pt\vrule
width2.5pt
height.25pt depth.15pt\vrule width.4pt
height6.65pt depth.15pt\ }

\def\proof{\noindent {\bf Proof. }}

\def\HollowBoxx #1#2#3{{\dimen0=#1 \advance\dimen0 by -#2       
       \dimen1=#1 \advance\dimen1 by #3                       
        \vrule height 0pt depth #3 width #2                   
       \hskip -#3
       \vrule height #1 depth #3 width #3}}                   
 \def\LeftContraction{\mathord{\kern1.45pt \HollowBoxx{6pt}{3.5pt}{.4pt}}\,}

 \def\HollowBox #1#2#3{{\dimen0=#1 \advance\dimen0 by -#3       
       \dimen1=#1 \advance\dimen1 by #3                       
        \vrule height #1 depth #3 width #3                    
        \vrule height 0pt depth #3 width #2                   
        \hskip -#3}}                                                   
 \def\RightContraction{\mathord{\, \HollowBox{6pt}{3.1pt}{.4pt}} \kern1.6pt}

\def\c{\LeftContraction\, }

\def \Omegabar{\overline \Omega}

\def \d{\partial}

\def\zbar{\overline{z}}
\def\dbar{\overline{\partial}}
\def \hbar{\overline{h}}

\begin{document}

\title{On  Analytic  Solvability and Hypoellipticity\\ 
For  $\dbar$ and $\dbar_b$}
\author{Michael Christ\thanks{Research 
supported in part by NSF grant DMS--9306833 and Macquarie University.} 
\and Song-Ying Li\thanks{Research supported in part by
NSF grant DMS--9500758.}}

\date{(November 11, 1995)}
\maketitle

\section{Introduction}

For each $0\le p,q\le n$, we denote by $\Lambda^{p,q}(\CC^n)$ 
the bundle of $(p,q)$-forms on $\CC^n$. Let $U$ be any open set in
$\CC^n$. We denote by ${\cal E}^{p,q}(U)$ the space of smooth sections of 
$\Lambda^{p,q}(\CC^n)$ over $U$. For any smooth hypersurface $M$
in $\CC^n$, we denote by $I^{p,q}(M)$ the ideal in 
$\Lambda^{p,q}(\CC^n)$ generated by $\rho$ and $\dbar \rho$,
where $\rho: \CC^n\to \RR$ is any smooth function that vanishes on $M$ and 
whose gradient is nowhere zero on $M$.
Let $\Lambda^{p,q}(M)$ denote
the orthogonal complement, with respect to the standard Hermitian metric
on $\CC^n$, of $I^{p,q}|_M$ in $\Lambda^{p,q}(\CC^n)|_M$. 
 
Let $\Omega$ be a bounded pseudoconvex domain in $\CC^n$ with
real analytic boundary $\d \Omega$, equipped with some
defining function $\rho\in C^\omega$ in a neighborhood of $\d\Omega$.
Let $L^2(\d \Omega)$ denote the Lebesgue space of measurable complex-valued 
functions on $\d \Omega$ that are square integrable with
respect to surface measure, and let ${\H}^2\subset L^2(\d\Omega)$ be
the subspace of all CR functions, that is, of functions annihilated
by $\dbar_b$. 
Let $S: L^2(\d \Omega)\to
{\H}^2$ be the Szeg\"o projection.  
Let $L^2(\Omega)$
be the Lebesgue space over $\Omega$, with respect to Lebesgue measure on $\CC^n$,
and let $A^2(\Omega)$ be the subspace of all holomorphic functions.
Let $P: L^2(\Omega) \to A^2(\Omega)$ be the Bergman projection. 
We denote by
$C^{\omega}(\d \Omega)$ the space of all real analytic functions over $\d \Omega$, and by
$C^{\omega} (\Omegabar)$ the space of all real analytic functions over $\Omegabar$. Let
$C^{\omega}_{(p,q)}(\Omegabar)$ denote the space of all $(p,q)$-forms with coefficients in
$C^{\omega}(\Omegabar)$, and let $C^{\omega}_{(p,q)}(\d \Omega)$ denote all forms in
$\Lambda^{p,q}(\d \Omega)$ with real analytic coefficients on $\d \Omega$.

There exist [CH1] bounded
pseudoconvex domains $\Omega\subset\CC^2$ with real
analytic boundaries, for which the Szeg\"o projection does not map 
$C^{\omega}(\d\Omega)$ to $C^{\omega}(\d\Omega)$.
Although it has not yet been proved, we think it highly likely that
the Bergman projection fails to preserve $C^\omega(\Omegabar)$ for those
same domains.
An analogous counterexample for the $C^\infty$ category has been
established more recently [CH2].

As is well known, the Szeg\"o and Bergman projections
are closely related to the $\dbar_b$ and 
$\Box_b$ 
equations and to the $\dbar$-Neumann
problem, respectively. Let $N_b$ denote the tangential Neumann 
operator and $N$ the Neumann operator. Thus 
$\Box_b N_b=I$ 
on $\Lambda^{p,q}(\d \Omega)$ and 
$\Box N=I$ 
on $\Lambda^{p,q}(\Omega)$. Then these are related to the 
Szeg\"o and Bergman projections respectively by
the following formulas of Kohn [FK]:
$$
S=I-\dbar_b^*N_b \dbar_b, \quad P=I-\dbar^* N \dbar. 
$$
Here $\dbar^*_b$, $\dbar^*$ denote the operators adjoint to $\dbar_b$ 
and $\dbar$, respectively.\footnote{In the case of $\CC^2$, we define the
tangential Neumann operator by $N_b f=u$ where $\dbar_b\dbar_b^*u=\pi f$,
$u$ is orthogonal to the nullspace of $\dbar_b^*$, and
$\pi f$ is the orthogonal projection of $f$ onto the range of $\dbar_b$.
Thus $\dbar_b\dbar_b^*N_b=\pi$, rather than the identity operator.} 

The counterexample of [CH1] for the Szeg\"o projection thus implies that
for bounded, pseudoconvex, real analytic domains $\Omega\subset\CC^2$,
for a general $(0,1)$ form $f$ in $C^\omega(\d\Omega)$ belonging to the 
range of $\dbar_b$,
the canonical solution $u=\dbar_b^*N_b f$ of $\dbar_b u=f$
need not always belong to $C^\omega(\d\Omega)$.
It is then a natural question, posed to us by N.~Sibony, whether there exists
a better solution of $\dbar_b u=f$, that does belong to $C^\omega(\d\Omega)$.
This is only an issue in the context of global existence, since the
Cauchy-Kowalevska theorem guarantees local $C^\omega$ solvability
given $C^\omega$ data.
Since it is likely that analogous counterexamples to global real analytic
regularity will eventually be established in higher dimensions
and for the $\dbar$--Neumann problem, the corresponding questions
for those contexts are of interest as well.
One aim of this paper is to answer these questions in the affirmative.

\begin{theorem}  Let $\Omega$ be a bounded pseudoconvex 
domain in $\CC^n$ with real analytic boundary, and suppose
that $0\le p\le n$ and $0<q\le n$. Let 
$f\in C^{\omega}_{(p,q)}(\Omegabar)$ be $\dbar$-closed. Then 
there is $u\in C^{\omega}_{(p,q-1)}(\Omegabar)$ 
such that $\dbar u=f$ in $\Omegabar$.
\end{theorem}

\begin{theorem}  Let $\Omega$ be a bounded pseudoconvex 
domain in $\CC^n$ with real analytic boundary. Then 

(i) If $n\ge 3$ and $f\in C^{\omega}_{(0,1)}(\d \Omega)$ is $\dbar_b$-closed,
there exists $u\in C^{\omega}(\d \Omega)$ satisfying $\dbar_b u=f$.

(ii) If $n=2$ and  $f\in C^{\omega}_{(0,1)}(\d \Omega) \cap 
\hbox{Range\,}(\dbar_b)$,
then  there exists $u\in C^{\omega}(\d \Omega)$ satisfying $\dbar_b u=f$.
\end{theorem}
 
Our second aim is to examine the connection between global regularity for
the Szeg\"o and Bergman projections on the one hand, and
the corresponding Neumann operators on the other. In the $C^\infty$ category
there is already known rigorously to be a close connection between the two [BS].
Theorems 1.1 and 1.2 have the following direct consequences.
Denote by $N_{p,q}$ the Neumann operator 
and by $P_{p,q}$ the Bergman projection
acting on $(p,q)$ forms on $\Omega$.\footnote{$P_{p,q}$ 
is the orthogonal projection onto the nullspace of $\dbar$.}
Denote by $N_b$ the tangential Neumann operator acting
on $(0,1)$ forms on $\d\Omega$,
and by $S$ the Szeg\"o projection, acting on functions.

\begin{corollary} 
Let $\Omega$ be a bounded pseudoconvex 
domain in $\CC^n$ with real analytic boundary. Then 

(i) Suppose that $n\ge 2$, $0\le p\le n$, and $0<q\le n$.
Then $P_{p,q-1}$ preserves $C^{\omega}_{(p,q-1)}(\Omegabar)$
if and only if $\dbar^*\circ N_{p,q}$
maps $C^{\omega}_{(p,q)}(\Omegabar)\cap \hbox{Kernel\,}(\dbar)$
to $C^{\omega}_{(p,q)}(\Omegabar)$.

(ii) Suppose that $n=2$. 
Then $S$ preserves $C^{\omega}(\d \Omega)$ if and only
if  $\dbar_b^*\circ N_b$ maps $C^{\omega}_{(0,1)}(\d \Omega)
\cap \hbox{Range\,} (\dbar_b)$ to $C^{\omega}(\d \Omega)$.

(iii) Suppose that $n > 2$. 
Then $S$ preserves $C^{\omega}(\d \Omega)$ if and only
if  $\dbar_b^*\circ N_b$ maps $C^{\omega}_{(0,1)}(\d \Omega)
\cap \hbox{Kernel\,}(\dbar_b)$
to $C^{\omega}(\d \Omega)$.
\end{corollary}

The first author thanks the Mathematics department of Macquarie University
for its hospitality. 
This work was carried out in part during visits of the second author
to UCLA and UCI, whose Mathematics departments he thanks for their
support.

\section{Proof of Theorem 1.1}

In this section we shall prove Theorem 1.1. We require the following
theorem of Diederich and Fornaess [DF1] and [DF2].

\begin{lemma}  Let $\Omega$ be a 
pseudoconvex bounded domain in $\CC^n$ with real analytic
boundary. Then there exists a decreasing sequence $\{\Omega^j\}$ 
of strictly pseudoconvex domains with smooth 
boundaries $\d \Omega^j$ such that $\Omegabar=\cap_{j=1}^{\infty} \Omega^j$.
\end{lemma}

\noindent{\bf Proof of Theorem 1.1.}
Let $f=\sum_{I,J} f_{IJ}dz^I\wedge d\bar z^J$
be a $\dbar$-closed $(p,q)$ form in $C^\omega(\overline{\Omega})$,
where $0<q\le n$.
Then there exists $j$ such that each $f_{IJ}$ can be extended
to a real analytic function on $\overline{\Omega_j}$.
Denoting the extension also by $f_{IJ}$, the form
$f=\sum_{I,J} f_{IJ}dz^I\wedge d\bar z^J$
is still $\dbar$--closed, by analytic continuation.

Since $\Omega^j$ is a smoothly bounded strictly
pseudoconvex domain in $\CC^n$, the
$\dbar$--Neumann problem is solvable on $\Omega^j$. Thus
[FK] there exists
$g\in C^{\infty}_{(p,q-1)}(\overline{\Omega^j})$
satisfying $\Box g=f$ on $\Omega^j$, where $\Box$
equals the ordinary Laplacian, acting componentwise on $g$.
$g$ also satisfies the $\dbar$--Neumann boundary conditions.

Since $\Box g\in C^\omega$, it follows (from ellipticity
of $\Box $, or from the explicit expression for the
fundamental solution for the Laplacian) that
$g\in C^\omega$ on the interior of $\Omega^j$.
In particular, $g\in C^\omega(\overline{\Omega})$,
and thus $\dbar^* g$ is also real analytic on
$\overline{\Omega}$.
Since $\dbar f=0$, it follows from the $\dbar$--Neumann
formalism [FK] that $\dbar(\dbar^*g)=f$.
\epf 

\section{Bochner extension for forms}

The strategy for the proof of Theorem 1.2, given a form
defined on $\d\Omega$, will be to extend it to a $\dbar$--closed 
form defined in a neighborhood of $\d\Omega$, or perhaps even
in a neighborhood of  
$\Omegabar$, then to argue as before. 

\begin{proposition}  Let $n\ge 2$ and $0\le q<n$. 
Let $\Omega\subset\CC^n$ be a bounded domain
with real analytic boundary. Let 
$f\in C^{\omega}_{(0,q)}(\d \Omega)$ be $\dbar_b$-closed. 
Then there exists a $\dbar$-closed form $F \in C^\omega_{(0,q)}$, 
defined in some neighborhood of $\d \Omega$,
whose restriction to $\d\Omega$, in the sense of forms, equals $f$.
\end{proposition}

Some related results for the $C^\infty$ category appear
in [FK] and [S].
Note that for $\CC^2$, the hypothesis that $\dbar_b f=0$ holds automatically
for all $(0,1)$ forms $f$ on $\d\Omega$.

Let $\rho$ be a defining function for
$\Omega$ that is real analytic in a neighborhood of $\d\Omega$. 
In particular, 
$\Omega=\{z:\rho(z)<0\}$.  
Let $N={4\over |\nabla \rho|^2} \sum_{j=1}^n
{\d \rho \over \d z_j} {\d \over \d \zbar_j}$. Let $\c$ be
the contraction operator on differential forms. The proof of
the next lemma is straightforward and is left to the reader.

\begin{lemma}  The following identities
hold in a neighborhood of $\d\Omega$ for any form $\phi$.
\vskip 3pt

(i) $\phi=N{\c} (\dbar \rho \wedge \phi)+\dbar \rho \wedge (N{\c} \phi)$

(ii) $N{\c} (N{\c}\phi) =0$.

\end{lemma} 

\noindent{\bf Proof of Proposition 3.1}.

To simplify notation we treat only the case $q=1$ explicitly,
but it will be apparent that the same reasoning applies to forms
of arbitrary degree.
Let a $\dbar_b$--closed $(0,1)$ form $f\in C^\omega(\d\Omega)$ be given.
We claim that in a sufficiently
small neighborhood of $\d\Omega$ there exists
a unique $C^\omega$ form $F$ satisfying
the differential equation 
$$ N\c\dbar F=0    \leqno(3.1)$$
with the side condition
$$N\c F=0,     \leqno(3.2)$$
whose restriction to $\d\Omega$ equals $f$.
Indeed, locally near any point in $\d\Omega$ we may fix $C^\omega$
$(0,1)$ forms $\{\bar\omega_j: 1\le j\le n-1\}$ that constitute an orthonormal
basis for the orthocomplement of $\dbar\rho$. The condition
$N\c F=0$ means that $F=\sum_{j\le n-1} F_j \bar\omega_j$ for
some coefficients $F_j$. Then $N\c\dbar F$ belongs also to the
span of the $\bar\omega_j$ and has the form 
$$
N\c\dbar F
= \sum_{j=1}^{n-1} \pm (\bar L F_j) \bar\omega_j
\leqno(3.3)
$$
modulo an operator of order $0$ acting on $F$, where $\bar L$ is
a complex vector field that is nowhere tangent to $\d\Omega$.
Thus (3.3) is a first-order system for which $\d\Omega$ is a
noncharacteristic hypersurface.
The coefficients $F_j$ are determined on $\d\Omega$ by $f$.
Therefore the Cauchy-Kowalevska theorem guarantees local existence
and uniqueness of a solution $F$ of (3.1) and (3.2) that equals $f$
on $\d\Omega$. Existence in a neighborhood of the full boundary 
follows from the local existence together with local uniqueness.

This form $F$ is $\dbar$-closed. Indeed,
$\dbar F = N\c(\dbar\rho\wedge\dbar F) + \dbar\rho\wedge (N\c\dbar F)$, 
and the second term vanishes by (3.1). On the other hand,
$G=N\c\dbar\rho\wedge \dbar F$ satisfies
the differential equation
$$
N\c\dbar G = N\c\dbar \, \dbar F = 0.
$$
But $G$ may be expressed as
$\sum_{i<j\le n-1} G_{ij} \bar\omega_i\wedge\bar\omega_j$
for certain $C^\omega$ coefficients $G_{ij}$. Just as above,
$$
N\c\dbar G = \sum_{i<j\le n-1}
\pm (\bar L G_{ij}) \bar\omega_i\wedge\bar\omega_j
$$
modulo a term of order zero belonging to the span of
the $\bar\omega_i\wedge\bar\omega_j$, where $\bar L$
is the same vector field as above. 
On $\d\Omega$, $G=\dbar_b f$ vanishes, so the uniqueness conclusion
of the Cauchy-Kowalevska theorem guarantees that 
$G=N\c\dbar\rho\wedge \dbar F$ vanishes identically in a neighborhood
of $\d\Omega$.
\epf 

\section{Proof of Theorem 1.2}
We will need the following well known existence theorem.

\begin{lemma} Suppose that $n\ge 3$ and that
$\Omega_j$ ($j=1,2$) are bounded, strictly
pseudoconvex domains in $\CC^n$  with smooth boundaries,
with $\Omega_1\subset \subset \Omega_2$.
Then for each $\dbar$--closed
$(0,1)$ form $f\in C^\infty(\overline{\Omega_2\setminus \Omegabar_1})$ 
there exists
a function $u\in C^\infty(\overline{\Omega_2\setminus \Omegabar_1})$ 
satisfying $\dbar u=f$.
\end{lemma} 

\proof Since $n\ge 3$, the Levi form has 
$n-1$ positive eigenvalues at each point of the outer boundary,
and has $n-1\ge 2$ negative eigenvalues at each point of the inner boundary. 
Thus $\Omega_2\setminus \Omegabar_1$
satisfies the $Z(1)$ condition at each boundary point,
and this condition is sufficient [FK] to guarantee the existence of a
smooth solution. 
\epf

\medskip
It is now straightforward to deduce Theorem 1.2 in the
easier case $n\ge 3$.
\medskip
 
\proof Let $f\in C^{\omega}_{(0,1)}(\d \Omega)$ satisfy $\dbar_b f=0$. 
By Lemma 2.1 and Proposition 3.1,
there are strictly pseudoconvex domains $\Omega_1$ and $\Omega_2$ such that
$\Omega_1\subset \subset \Omega \subset \subset \Omega_2$ and a $\dbar$-closed
$(0,1)$-form  $\phi=\phi_1 d\zbar_1+\cdots + \phi_n d\zbar_n\in C^{\omega}(\bar U)$, 
where $U= (\Omega_2\setminus \Omegabar_1)$, such that
$$
N{\c}(\dbar \rho\wedge  \phi ) =f\quad \hbox{ on }\  
\d \Omega.
$$
By Lemma 4.1, there is $u \in C^{\infty}(\Omega_2\setminus \Omega_1)$
such that $\dbar u=\phi$. The ellipticity of $\dbar$ implies
that $u \in C^{\omega}(\Omega_2\setminus \Omega_1)$, and 
$$
\dbar_b u=N{\c} (\dbar \rho \wedge
\dbar u)=N{\c}(\dbar \rho \wedge \phi)=f\quad 
\hbox{ on  }\d \Omega.
$$
The proof of Theorem 1.2 is thus complete when $n\ge 3$.
\epf\medskip

Next we shall prove Theorem 1.2 for the case $n=2$.

\begin{proposition} Let $\Omega$ be a bounded domain in $\CC^2$ with
real analytic boundary. Let $F$ be a 
$\dbar$-closed $(0,1)$-form with $C^\infty$ coefficients, defined
in a neighborhood of $\d\Omega$, whose restriction to $\d\Omega$
belongs to $\dbar_b(C^\infty)$. 
Then there exists a $\dbar$-closed $(0,1)$-form $\tilde{F}$ having
Lipschitz continuous coefficients, defined on 
a neighborhood $\Omega'$ of $\Omegabar$, that agrees with $F$ on
$\Omega'\setminus \Omegabar$. 
\end{proposition}

\proof 
Fix a $C^\omega$ defining function $\rho$ for $\Omega$, 
and fix $u_0\in C^\infty(\d\Omega)$ satisfying $\dbar_b u_0
=N\c\dbar\rho\wedge F$ on $\d\Omega$.
Extend $u_0$  to a smooth function on $\Omega'$.  
Then
$$
\dbar u_0 \wedge \dbar \rho - F\wedge \dbar \rho
=O(\rho).\leqno(4.1)
$$
Let 
$$
h=N{\c} F-N{\c} (\dbar u_0),
$$
and define
$$
u(z)= u_0(z)+\rho(z) h(z).  \leqno(4.2)
$$ 
Then
$$
\dbar u=\dbar u_0(z)+ h(z) \dbar \rho(z)+ \rho(z) \dbar h(z)
$$
It is easy to see that
$$
\dbar u \wedge \dbar \rho =F\wedge \dbar \rho\quad \hbox{ on } \d \Omega,
\leqno(4.3)
$$
and
$$
N{\c} (\dbar u)-N{\c} F=\rho(z) N{\c}(\dbar h)=
\rho(z) (N{\c} \dbar) (N{\c}(F-\dbar u_0)).\leqno(4.4)
$$
We let 
$$
\tilde{F}=\cases{ F & \quad on  $\Omega'\setminus\Omega$ \cr
\dbar u & \quad on  $\Omega$. \cr}
$$
By (4.1) and (4.4), $\dbar u $ and $F$ agree to first order
$\d \Omega$. Thus $\tilde F$ has Lipschitz coefficients. 
$\dbar \tilde{F}\in L^\infty$ vanishes identically on both
$\Omega$ and $\Omega' \setminus \Omegabar$, hence vanishes on $\Omega'$
in the sense of distributions. The
proof of the proposition is complete.
\epf

\medskip
To complete the proof of Theorem 1.2, we need the following lemma.

\begin{lemma}
Let $\Omega\subset\CC^2$ be bounded and pseudoconvex, with $C^\omega$
boundary. Suppose that $\Omega'$ is a neighborhood of $\Omegabar$,
that $F$ is a $\dbar$--closed $(0,1)$
form having $C^\omega$ coefficients in a neighborhood of $\d\Omega$,
and that $v\in C^1(\Omega'\setminus\Omega$)
satisfies $\dbar v=F$.
Then $v$ extends to a function belonging to $C^\omega$ in a neighborhood
of $\d\Omega$.
\end{lemma}

\proof
$v\in C^\omega(\Omega'\setminus\Omegabar)$, since $\dbar v\in C^\omega$
there. 
For any point $p\in \d \Omega$
there exist $\delta>0$  and a real analytic function $u_p$ satisfying
$\dbar u_p =F$ on the ball $B(p,\delta)$. 
Thus $\dbar (v-u_p)=0$ on $B(p,\delta) \cap  (\Omega' \setminus \Omega) $.
Therefore $v-u_p$ is holomorphic in 
$B(p,\delta)\cap(\Omega'\setminus \Omegabar)$.

$p$ is a point of finite type in $\d\Omega$, since
$\Omega\subset\CC^2$ is bounded and has real analytic boundary.
$\dbar_b(v-u_p)=0$ on $\d\Omega$ near $p$, by continuity
since $v$ was assumed to be $C^1$ up to the boundary.
Since $\Omega$
is pseudoconvex, and $p$ is a boundary point of finite type,
a theorem of Bedford and Fornaess [BF] asserts that for some small
ball $B$ centered at $p$ there
exists a holomorphic function $h$ on $B\cap\Omega$, continuous
in $B\cap\Omegabar$, that agrees with $v-u_p$ on $\d\Omega$.
Extend $h$ to $B$ by defining it to be $v-u_p$ in $B\setminus
\Omegabar$.

Now $h$ is continuous, and is annihilated by $\dbar$ both inside
and outside $\d\Omega$, near $p$. Therefore $\dbar h=0$ in
the sense of distributions, so $h$ is a genuine holomorphic
function near $p$. It agrees with $v-u_p$ in $B\setminus\Omegabar$,
hence $v-u_p$ has been extended to a function holomorphic in
a neighborhood of $p$.
Since $u_p$ is real analytic, $v$ itself has thus been extended to a
function
real analytic in a neighborhood of $p$, which was an arbitrary point of
$\d\Omega$.
\epf 

\medskip
We are now in a position to complete the proof of Theorem 1.2 for the case
$n=2$.

\proof Let $f\in C^{\omega}_{(0,1)}(\d \Omega) \cap \hbox{Range}(\dbar_b)$.
By Proposition 3.1, there exists a $\dbar$--closed $(0,1)$
form $F$, defined in a neighborhood of $\d\Omega$, that extends
$f$ in the sense of forms.
By Proposition
4.2, there are a neighborhood $\Omega'$ of $\Omegabar$ and a
$\dbar$-closed $(0,1)$-form 
$\tilde{F}$ 
with Lipschitz continuous coefficients such that
$\tilde{F}=F$ on $\Omega'\setminus \Omega$. 
By Lemma 2.1, we may assume $\Omega'$ to be strictly pseudoconvex.
Therefore there exists a  solution
$v$ of $\dbar v=\tilde F$ in $\Omega'$,
and elliptic regularity theorems guarantee that
$v\in C^{1,\alpha}$, for every $\alpha<1$.
In particular, $\dbar v=F$ in $\Omega'\setminus\Omegabar$. 
By Lemma 4.3, $v\in C^{\omega}(\d\Omega)$ 
and consequently $\dbar v\equiv F$ in a neighborhood of
$\d\Omega$, by analytic continuation.
Since $F$ extends $f$, we thus have $\dbar_b v=f$ on $\d \Omega$,
concluding the proof of Theorem 1.2 for the case $n=2$.
\epf

\section{Proof of Corollary 1.3}

\proof 
Write $P=P_{p,q-1}$, $N=N_{p,q}$.
>From the formula $P= I -\dbar^* N \dbar$ it follows
directly that if $\dbar^*\circ N$ maps
$C^{\omega}(\Omegabar)\cap \hbox{ker}(\dbar)$ to
$C^{\omega}(\Omegabar)$, then $P$ preserves
$C^{\omega}(\Omegabar)$. 

Conversely,
for any $\dbar$-closed $(p,q)$-form $f\in C^{\omega}(\Omegabar)$,  
Theorem 1.1 guarantees the existence of 
$u\in C^{\omega}_{(p,q-1)}(\Omegabar)$ satisfying
$ \dbar u=f$.
Then $Pu\in C^{\omega}(\Omegabar)$, by the global regularity
hypothesis on $P$. Since $u-Pu$ is orthogonal to 
the $L^2$ nullspace of $\dbar$, it equals $\dbar^* Nf$.

(ii) and (iii) can be proved similarly by using Theorem 1.2  
and the formula $S = I-\dbar_b^*\circ N_b\circ\dbar_b$.
Existence of $N_b$ and subelliptic estimates follow from
work of Kohn [K], 
since boundaries of bounded, pseudoconvex, real analytic
domains are necessarily of finite ideal type as defined in [K].
\epf

\bigskip

\noindent Department of Mathematics, UCLA,
Los Angeles, CA 90095-1555.

\noindent e-mail: christ@math.ucla.edu

\medskip

\noindent Department of Mathematics, Washington University, 
St. Louis, MO 63130. 

\noindent e-mail: songying@math.wustl.edu
\end{document}